\documentclass[12pt]{article}

\newtheorem{definition}{Definition}[section]
\newtheorem{theorem}[definition]{Theorem}
\newtheorem{lemma}[definition]{Lemma}
\newtheorem{corollary}[definition]{Corollary}
\newtheorem{remark}[definition]{Remark}
\newtheorem{example}[definition]{Example}
\newtheorem{conjecture}[definition]{Conjecture}
\newtheorem{problem}[definition]{Problem}
\newtheorem{note}[definition]{Note}

\usepackage{amssymb}
\usepackage{latexsym}

\def\C{\mathbb C}

\newcommand{\fld}{{\mathbb K}}
\newcommand{\K}{{\mathbb K}}

\newcommand{\beast}{\begin{eqnarray*}}
\newcommand{\eeast}{\end{eqnarray*}}

\begin{document}

\title{ \bf The shape of    
 a tridiagonal pair}
\author{Tatsuro Ito and Paul Terwilliger  
}
\date{}
\maketitle
\begin{abstract} Let $\fld $ denote an algebraically closed field
with characteristic 0. 
Let $V$ denote a vector space
over  $\fld$ with finite positive dimension
and let $A,A^*$ denote a tridiagonal pair on $V$.
We make an assumption about this pair.
Let $q$ denote a nonzero scalar in $\K$ which is
not a root of unity. 
We assume $A$ and $A^*$ satisfy the $q$-Serre relations
\begin{eqnarray*}
A^3A^* - \lbrack 3 \rbrack A^2A^*A
+ \lbrack 3 \rbrack AA^*A^2 - A^*A^3&=&0, 
\\
A^{*3}A - \lbrack 3 \rbrack A^{*2}AA^*
+ \lbrack 3 \rbrack A^*AA^{*2} - AA^{*3}&=&0, 
\end{eqnarray*}
where $\lbrack 3\rbrack=(q^3-q^{-3})/(q-q^{-1})$.
Let $(\rho_0, \rho_1,\ldots,\rho_d)$ denote the shape vector
for $A,A^*$. We show the entries in this shape vector are bounded
above by binomial coefficients as follows:
\begin{eqnarray*}
\rho_i \;\leq \;\biggl( {{d}\atop {i}} \biggr) \qquad \qquad (0 \leq i \leq d).
\end{eqnarray*}
We obtain this result by displaying a spanning set for $V$.
\end{abstract}


\section{Introduction}

\noindent 
Throughout this paper,
 $\K$ will denote a field and  $V$ will denote
a vector space over $\K$ with finite positive dimension.
 Let $A:V\rightarrow V$ denote a linear transformation and
 let $W$ denote a subspace of $V$. We call $W$ an {\it eigenspace} of $A$ whenever $W\not=0$ and there exists $\theta \in \K$ such that 
\begin{eqnarray*}
W=\lbrace v \in V \;\vert \;Av = \theta v\rbrace.
\end{eqnarray*}
We say $A$ is {\it diagonalizable} whenever
$V$ is spanned by the eigenspaces of $A$.

\medskip
\noindent We now recall the notion of a {\it tridiagonal pair}.

\begin{definition}  
\cite{TD00}
\label{def:tdp}
\rm
By a {\it tridiagonal pair}  on $V$,
we mean an ordered pair $A,A^*$ where
$A:V\rightarrow V$ and 
$A^*:V\rightarrow V$ 
 are linear transformations which satisfy
the following four conditions.
\begin{enumerate}
\item Each of $A,A^*$ is diagonalizable.
\item There exists an ordering $V_0, V_1,\ldots, V_d$ of the  
eigenspaces of $A$ such that 
\begin{equation}
A^* V_i \subseteq V_{i-1} + V_i+ V_{i+1} \qquad \qquad (0 \leq i \leq d),
\label{eq:t1}
\end{equation}
where $V_{-1} = 0$, $V_{d+1}= 0$.
\item There exists an ordering $V^*_0, V^*_1,\ldots, V^*_\delta$ of
the  
eigenspaces of $A^*$ such that 
\begin{equation}
A V^*_i \subseteq V^*_{i-1} + V^*_i+ V^*_{i+1} \qquad \qquad (0 \leq i \leq \delta),
\label{eq:t2}
\end{equation}
where $V^*_{-1} = 0$, $V^*_{\delta+1}= 0$.
\item There does not exist a subspace $W$ of $V$ such  that $AW\subseteq W$,
$A^*W\subseteq W$, $W\not=0$, $W\not=V$.
\end{enumerate}
\end{definition}

\begin{note} \rm
According to a common 
notational convention, $A^*$ denotes the conjugate-transpose
of $A$.
 We are not using this convention.
In a tridiagonal pair $A,A^*$ the linear transformations 
$A$ and $A^*$ are arbitrary subject to (i)--(iv) above.
\end{note}

\noindent
In order to motivate our results we recall a few facts about 
tridiagonal pairs.
Let $A,A^*$ denote a tridiagonal pair on $V$
and let the integers $d, \delta$ be as in
Definition \ref{def:tdp}(ii), (iii) respectively.  By 
 \cite[Lemma 4.5]{TD00}
 we have
$d=\delta$; we call this common value the {\it diameter}
of $A,A^*$.  
An ordering of the eigenspaces of $A$ (resp. $A^*$)
is called
{\it standard} whenever it satisfies
(\ref{eq:t1}) (resp. 
(\ref{eq:t2})). 
We comment on the uniqueness of the standard ordering.
Let 
 $V_0, V_1,\ldots, V_d$ denote a standard ordering
of the eigenspaces of $A$.
Then the ordering 
 $V_d, V_{d-1},\ldots, V_0$ 
is standard and no other ordering is standard.
A similar result holds for the 
 eigenspaces of $A^*$.
Let
 $V_0, V_1,\ldots, V_d$
(resp.  $V^*_0, V^*_1,\ldots, V^*_d$)
denote a standard ordering
of the eigenspaces of $A$ (resp. $A^*$).
By \cite[Corollary 5.7]{TD00},  for $0 \leq i \leq d$ 
the spaces $V_i$,  $V^*_i$ have the same dimension;
we denote this common dimension by $\rho_i$.
By the construction $\rho_i\not=0$.
 By  \cite[Corollary 5.7]{TD00} 
and \cite[Corollary 6.6]{TD00},
the sequence $\rho_0, \rho_1, \ldots, \rho_d$
is symmetric and unimodal; that is
$ \rho_i =\rho_{d-i}$ for $0 \leq i \leq d$
and 
$\rho_{i-1} \leq \rho_{i}$ for  $1 \leq i \leq d/2$.
We refer to the sequence
 $(\rho_0, \rho_1, \ldots, \rho_d)$ as the {\it shape vector} of $A,A^*$.

\medskip
\noindent
The following special case has received a lot of attention.
By a {\it Leonard pair} we mean a tridiagonal pair
which has shape vector $(1,1,\ldots, 1)$.
There is a natural correspondence between
the Leonard pairs and a family of orthogonal polynomials
consisting of
the $q$-Racah polynomials
 \cite{AWil}
and some related polynomials
of the Askey-scheme 
\cite{KoeSwa}, 
\cite{TLT:array},
\cite{TLTqrac}.
There is a classification of
Leonard pairs in
\cite{LS99}. This classification amounts to
a linear algebraic version of a theorem of
D. Leonard
\cite{BanIto}, \cite{Leon}
concerning the $q$-Racah polynomials.
See 
\cite{TersubI},
\cite{qSerre},
\cite{LS24},
\cite{conform},
\cite{lsint},
\cite{Terint},
\cite{TLT:split},
\cite{TLT:array},
\cite{TLTqrac}
for more information about Leonard pairs.

\medskip
\noindent In this paper our focus is on general tridiagonal pairs.
We will discuss the following conjecture.

\begin{conjecture} (Ito, Tanabe, Terwilliger \cite{TD00})
 \label{con:bc}
 Let  $(\rho_0, \rho_1, \ldots, \rho_d)$ denote the shape vector
of a tridiagonal pair.
Then
\begin{eqnarray*}
\rho_i \;\leq \;\biggl( {{d}\atop {i}} \biggr) \qquad \qquad (0 \leq i \leq d).
\end{eqnarray*}
\end{conjecture} 

\noindent In the present paper we will prove Conjecture
 \ref{con:bc} for a certain type of tridiagonal pair. We will describe this
type shortly. For now we review some more facts about general 
tridiagonal pairs.

\medskip
\noindent
Let $A,A^*$ denote a tridiagonal pair on $V$.
By
\cite[Theorem 10.1]{TD00}
there exists a sequence
$\beta, \gamma, \gamma^*, \varrho,
 \varrho^*$ of scalars taken
from $\K$ such that 
both
\begin{eqnarray}
\lbrack A,A^2A^*-\beta AA^*A + 
A^*A^2 -\gamma (AA^*+A^*A)-\varrho A^*\rbrack &=&0, 
\label{eq:m1n}
\\
\lbrack A^*,A^{*2}A-\beta A^*AA^* + AA^{*2} -\gamma^* (A^*A+AA^*)-
\varrho^* A\rbrack&=&0,
\label{eq:m2n}
\end{eqnarray}
where $\lbrack r,s \rbrack $ means $\,rs-sr$.
The sequence is unique if the diameter is at least $3$.
The equations 
(\ref{eq:m1n}),
(\ref{eq:m2n}) 
are known as the {\it tridiagonal relations}
\cite{TD00}, \cite{qSerre}.
An ordering of the eigenvalues of $A$ (resp. $A^*$)
will be called {\it standard}
whenever the corresponding ordering of the eigenspaces of $A$
(resp. $A^*$) is standard.
Let $\theta_0,\theta_1,\ldots, \theta_d$ (resp. 
 $\theta^*_0,\theta^*_1,\ldots, \theta^*_d$)
 denote a standard ordering
of the eigenvalues of $A$ (resp. $A^*$). By
\cite[Theorem 11.1]{TD00} we have
\begin{eqnarray}
\theta^2_{i-1}
-
\beta \theta_{i-1}\theta_i
+
\theta^2_i
-\gamma (\theta_{i-1}+\theta_i)
&=&\varrho 
\qquad \qquad 
(1 \leq i \leq d),
\label{eq:eigrec}
\\
\theta^{*2}_{i-1}
-
\beta \theta^*_{i-1}\theta^*_i
+
\theta^{*2}_i
-\gamma^* (\theta^*_{i-1}+\theta^*_i)
&=&\varrho^*
\qquad \qquad 
(1 \leq i \leq d).
\label{eq:eigrec2}
\end{eqnarray}
\noindent 
 By \cite[Theorem 4.6]{TD00}
there exists a unique sequence $U_0, U_1, \ldots, U_d$ consisting of
subspaces of $V$ such that
\begin{eqnarray}
&&\qquad V = U_0+ U_1 + \cdots + U_d \qquad \qquad (\hbox{direct sum}),
\label{eq:vdec}
\\
&&(A-\theta_iI)U_i\subseteq U_{i+1} \qquad \quad (0 \leq i <d),  
\qquad (A-\theta_dI)U_d=0,   
\label{eq:Aonui}
\\
&& (A^*-\theta^*_iI)U_i\subseteq U_{i-1} \qquad \quad (0 < i \leq d), 
\qquad (A^*-\theta^*_0I)U_0=0.    
\label{eq:Asonui}
\end{eqnarray}
By  \cite[Corollary 5.7]{TD00}   
the space $U_i$ has dimension $\rho_i$ for $0 \leq i \leq d$,
where $(\rho_0, \rho_1, \ldots, \rho_d)$ is the shape vector
for $A,A^*$.
We call the sequence $U_0,U_1,\ldots, U_d$
the {\it split decomposition} for $A,A^*$ (with respect to
the orderings $\theta_0,\theta_1,\ldots, \theta_d$ and
 $\theta^*_0,\theta^*_1,\ldots, \theta^*_d$).
For $0 \leq i \leq d$ let 
$F_i:V\rightarrow V$ denote the linear transformation which satisfies both
\begin{eqnarray}
&&\qquad (F_i-I)U_i=0,
\label{eq:fi1}
\\
&&F_iU_{j}=0 \quad \hbox{if}\quad j\not=i, \qquad (0 \leq j \leq d).
\label{eq:fi2}
\end{eqnarray}
In other words $F_i$  is the projection map from $V$ onto $U_i$.
We observe
\begin{eqnarray*}
&&F_iF_j = \delta_{ij}F_i \qquad \qquad (0 \leq i,j\leq d),
\label{eq:fifiidem}
\\
&&\qquad F_0+F_1+\cdots +F_d = I,
\label{eq:fisumtoident}
\\
&&F_iV = U_{i} \qquad \qquad (0 \leq i \leq d).
\label{eq:fivjustvi}
\end{eqnarray*}
In view of 
(\ref{eq:Aonui}), 
(\ref{eq:Asonui})
we define
\begin{eqnarray}
R&=&A-\sum_{h=0}^d \theta_hF_h,
\label{eq:defR}
\\
L&=&A^*-\sum_{h=0}^d \theta^*_hF_h
\label{eq:defL}
\end{eqnarray}
and observe
\begin{eqnarray}
&& RU_i \subseteq U_{i+1} \qquad \qquad (0 \leq i <d), \qquad RU_d=0,
\label{eq:raction}
\\
&&LU_i \subseteq U_{i-1} \qquad \qquad  (0 < i \leq d), \qquad LU_0=0.
\label{eq:laction}
\end{eqnarray}
Combining 
(\ref{eq:raction}),
(\ref{eq:laction}) with
(\ref{eq:vdec}) we find
\begin{eqnarray}
\label{eq:nil}
R^{d+1} = 0, \qquad \qquad 
L^{d+1} = 0.
\end{eqnarray}
We call $R$ (resp. $L$)
the {\it raising map} (resp. {\it lowering map}) for
$A,A^*$ 
with respect to $U_0, U_1, \ldots, U_d$.

\medskip
\noindent We now describe the type of 
tridiagonal pair for
 which we will prove
Conjecture \ref{con:bc}.

\begin{definition} 
\label{def:1p1}
\rm
For the rest of this paper, we assume $\K$ is algebraically closed
with characteristic 0.
We fix a nonzero scalar $q\in \fld$
which is not a root of unity.
We let $A,A^*$ denote a  tridiagonal pair on $V$
 which satisfies
(\ref{eq:m1n}),
(\ref{eq:m2n})
where
\begin{eqnarray}
\beta = q^2 + q^{-2},
\qquad \qquad  
\gamma = \gamma^* = \varrho= \varrho^* = 0. 
\label{eq:simpleassume}
\end{eqnarray}

\end{definition}
\begin{note}
\rm
The scalar $q$ which appears in \cite{TD00}   
is the same as $q^2$ in the present paper. We make this adjustment
for notational convenience.
\end{note}

\medskip
\noindent Before proceeding we discuss the significance of
assumption (\ref{eq:simpleassume}) from several points of view.
Let $A,A^*$ denote the tridiagonal pair
in Definition
\ref{def:1p1}.
Evaluating
 (\ref{eq:m1n}),
(\ref{eq:m2n})
using
(\ref{eq:simpleassume})
we obtain
\begin{eqnarray}
\label{eq:qs1}
 A^3A^* - \lbrack 3 \rbrack A^2A^*A
+ \lbrack 3 \rbrack AA^*A^2 - A^*A^3&=&0,
 \\
 A^{*3}A - \lbrack 3 \rbrack A^{*2}AA^*
+ \lbrack 3 \rbrack A^*AA^{*2} - AA^{*3}&=&0.
\label{eq:qs2} 
\end{eqnarray}
We are using the notation
\begin{eqnarray}
 \lbrack n\rbrack=\frac{q^n-q^{-n}}{q-q^{-1}}\qquad \qquad n=0,1,\ldots
\label{eq:qi}
\end{eqnarray}
The equations 
(\ref{eq:qs1}),
(\ref{eq:qs2}) are known as the {\it $q$-Serre relations}, and
are among the defining relations for the quantum affine algebra
$U_q(\widehat{ sl_2})$. See 
\cite{cp} and \cite{lu} for more information on $U_q(\widehat{ sl_2})$. 

\medskip
\noindent 
Assumption 
 (\ref{eq:simpleassume}) has the following significance 
for the eigenvalues.
Let $A,A^*$ denote the tridiagonal pair 
 in
Definition 
\ref{def:1p1}.
Let $\theta_0, \theta_1,\ldots, \theta_d$
(resp.
 $\theta^*_0, \theta^*_1,\ldots, \theta^*_d$)
denote a standard ordering of the eigenvalues of $A$ (resp. $A^*$).
Evaluating
(\ref{eq:eigrec})
using 
$\beta = q^2+q^{-2}$, 
$\gamma=0$, 
$\varrho=0$
we find either $\theta_i=q^2\theta_{i-1}$ for
$1 \leq i \leq d$
or $\theta_i=q^{-2}\theta_{i-1}$ for
$1 \leq i \leq d$. Replacing 
$\theta_0, \theta_1,\ldots, \theta_d$
by 
$\theta_d, \theta_{d-1},\ldots, \theta_0$ if necessary
we may assume 
 $\theta_i=q^2\theta_{i-1}$ for
$1 \leq i \leq d$. In this case
there exists
$a \in \fld$ such that
\begin{eqnarray}
\label{eq:eig}
\theta_i = a q^{2i-d} \qquad  (0 \leq i \leq d).
\end{eqnarray}
Similarly
replacing 
$\theta^*_0, \theta^*_1,\ldots, \theta^*_d$
by 
$\theta^*_d, \theta^*_{d-1},\ldots, \theta^*_0$ if necessary,
 there exists
$a^* \in\K$ such that
\begin{eqnarray}
\label{eq:eigs}
\theta^*_i = a^* q^{d-2i} \qquad \qquad (0 \leq i \leq d).
\end{eqnarray}

\noindent Assumption 
(\ref{eq:simpleassume}) has the following significance
for the raising and lowering maps.
Let $A,A^*$ denote the tridiagonal pair
in Definition 
\ref{def:1p1}. Let
$U_0, U_1, \ldots, U_d$ denote the split
decomposition
with respect to
(\ref{eq:eig}), 
(\ref{eq:eigs})
and let $R$ (resp. $L$)
denote the corresponding raising (resp. lowering)
map.
By \cite[Theorem 12.2]{TD00} 
the maps $R, L$ satisfy the $q$-Serre relations
\begin{eqnarray}
R^3L - \lbrack 3 \rbrack R^2LR+ \lbrack 3 \rbrack RLR^2 - LR^3&=&0, 
\label{eq:ms1}
\\
L^3R - \lbrack 3 \rbrack L^2RL+ \lbrack 3 \rbrack LRL^2 - RL^3&=&0.
\label{eq:ms2}
\end{eqnarray}
We emphasize that each of the pairs $R,L$ and $A,A^*$ satisfy the
$q$-Serre relations.  
Each of $R, L$
 is nilpotent on $V$ and each of $A, A^*$
is diagonalizable on $V$.

\medskip
\noindent
The  following is our main result.
\begin{theorem} 
\label{thm:bc}
Let $A,A^*$ denote the tridiagonal pair in
Definition 
\ref{def:1p1} and let 
$(\rho_0,\rho_1,\ldots, \rho_d)$ denote the 
corresponding shape vector. Then
\begin{eqnarray*}
\rho_i \;\leq \;\biggl( {{d}\atop {i}} \biggr) 
\qquad \qquad (0 \leq i \leq d).
\end{eqnarray*}
\end{theorem} 

\medskip
\noindent 
In order to prove
Theorem \ref{thm:bc}, for $0\leq i \leq d$
we display a spanning set for $U_i$ consisting
of 
$({{d}\atop{i}})$ vectors. We proceed as follows.

\begin{theorem} 
\label{con:Vspanguess}
Let $A,A^*$ denote the tridiagonal pair in 
Definition 
\ref{def:1p1}.
Let $U_0, U_1, \ldots, U_d$ denote the split
decomposition with respect to
(\ref{eq:eig}), 
(\ref{eq:eigs})
and let $R, L$
denote the corresponding raising and lowering
maps.
Let $v$ denote a nonzero vector in $U_0$.
Then $V$ is spanned by
the vectors of the form
\beast
L^{i_1}R^{i_2}L^{i_3}R^{i_4} \cdots R^{i_n}v,
\eeast
where 
$i_1, i_2, \ldots, i_n$ ranges over all sequences
such that $n$ is a nonnegative even integer, and 
 $i_1, i_2, \ldots i_n$ are integers
satisfying $0\leq i_1<i_2<\cdots <i_n\leq d$.
\end{theorem}

\medskip
\noindent We prove Theorem
\ref{con:Vspanguess}  at the end of Section 2.
In order to illustrate 
Theorem
\ref{con:Vspanguess} 
we now restate it
in concrete terms for $d=3$.

\begin{corollary} With reference to
Theorem
\ref{con:Vspanguess}, for
$d=3$
the space
$V$ is spanned by the vectors 
\beast
v,\quad  Rv, \quad R^2v, \quad  R^3v, 
\quad LR^2v, \quad LR^3v,
\quad L^2R^3v, \quad  
RL^2R^3v.
\eeast
\end{corollary}

\medskip
\noindent 
Combining Theorem
\ref{con:Vspanguess}
with 
(\ref{eq:raction}),
(\ref{eq:laction}) 
we routinely obtain the following theorem.

\begin{theorem} 
\label{thm:us}
Let $A,A^*$ denote the tridiagonal pair in 
Definition 
\ref{def:1p1}.
Let $U_0, U_1, \ldots, U_d$ denote the split
decomposition with respect to
(\ref{eq:eig}), 
(\ref{eq:eigs})
and let $R, L$
denote the corresponding raising and lowering
maps.
Let $v$ denote a nonzero vector in $U_0$.
Then for $0 \leq i \leq d$ the space $U_i$ is spanned by
the vectors of the form
\beast
L^{i_1}R^{i_2}L^{i_3}R^{i_4} \cdots R^{i_n}v,
\eeast
where 
$i_1, i_2, \ldots, i_n$ ranges over all sequences
such that (i) $n$ is a nonnegative even integer;
(ii)  $i_1, i_2, \ldots i_n$ are integers
satisfying $0\leq i_1<i_2<\cdots <i_n\leq d$ and $i=\sum_{h=1}^n i_h(-1)^h$.
\end{theorem}

\begin{corollary}
With reference to
Theorem
\ref{thm:us},
for $d=3$ 
the following (i)--(iv) hold.
\begin{enumerate}
\item
The space $U_0$ is spanned by
\beast
v.
\eeast
\item The space $U_1$ is spanned by
\beast
Rv, \quad LR^2v,\quad L^2R^3v.
\eeast
\item The space $U_2$ is spanned by
\beast
\quad R^2v, 
 \quad LR^3v,
 \quad  
RL^2R^3v.
\eeast
\item The space $U_3$ is spanned by
\beast
R^3v.
\eeast
\end{enumerate}
\end{corollary}

\section{The algebra $\mathcal A$}

\noindent Our next goal is to prove 
Theorem
\ref{con:Vspanguess}.
In order to do this we 
undertake a careful investigation of the 
$q$-Serre relations.
We begin with a definition.

\begin{definition}
\label{def:a}
\rm
We let $\mathcal A$ denote the associative $\fld$-algebra with identity
generated by symbols $x,y$ subject to the $q$-Serre relations
\begin{eqnarray*}
x^3y-\lbrack 3 \rbrack x^2yx+\lbrack 3 \rbrack xyx^2-yx^3&=&0,
\\
y^3x-\lbrack 3 \rbrack y^2xy+\lbrack 3 \rbrack yxy^2-xy^3&=&0.
\end{eqnarray*}
\end{definition}

\begin{remark}
\rm
The algebra $\mathcal A$ in Definition 
\ref{def:a}
is often called the
{\it positive part of 
$U_q({\widehat {sl_2}})$}.
\end{remark}

\begin{definition}
\label{def:s}
\rm
Referring to Definition \ref{def:a},
we let $\sigma_0:{\mathcal A}\rightarrow {\mathcal A}$ denote the $\fld$-algebra isomorphism
which sends $x$ to $y$ and $y$ to $x$.
We let $\sigma_1:{\mathcal A}\rightarrow{\mathcal A}$ denote the $\fld$-algebra antiisomorphism
which stabilizes each of 
$x,y$. The action of $\sigma_1$ inverts the word order.
We observe $\sigma_0, \sigma_1$  commute.
We define $\sigma=\sigma_0\sigma_1$ and
observe each of $\sigma_0, \sigma_1,\sigma$ is an involution.
\end{definition} 

\begin{example}
\label{ex:sig01}
Applying each of $\sigma_0, \sigma_1, \sigma$ to
$xyx^2y^2$ we get
\begin{eqnarray*}
yxy^2x^2,
\qquad  y^2x^2yx,
 \qquad x^2y^2xy
\end{eqnarray*}
 respectively.
\end{example}

\medskip
\noindent
We now state a lemma which we will find useful.
We first recall some notation.
Let $n$ denote a nonnegative integer.
With reference to 
(\ref{eq:qi})
we define
\begin{eqnarray*}
\lbrack n\rbrack !=
\lbrack n\rbrack 
\lbrack n-1\rbrack 
\cdots
\lbrack 2 \rbrack
\lbrack 1\rbrack.
\end{eqnarray*}
We interpret $\lbrack 0\rbrack !=1$.
We also define
\begin{eqnarray*}
\biggl\lbrack
{{n}\atop {i}}\biggr\rbrack
=\frac{\lbrack n\rbrack!}{\lbrack i\rbrack ! \lbrack n-i \rbrack !}
\qquad (0 \leq i\leq n).
\end{eqnarray*}

\begin{lemma}
\cite[Proposition 7.1.5]{lu}
\label{lem:ho}
The following relations
hold in $\mathcal A$.
Let $r$ denote a positive integer.
Then
\begin{eqnarray}
\sum_{i=0}^{2r+1}
(-1)^i
\biggl\lbrack
{{2r+1}\atop {i}}\biggr\rbrack
x^iy^rx^{2r+1-i} &=&0,
\label{eq:hoq}
\\
\sum_{i=0}^{2r+1}
(-1)^i
\biggl\lbrack
{{2r+1}\atop {i}}\biggr\rbrack
y^ix^ry^{2r+1-i} &=&0.
\label{eq:hoq2}
\end{eqnarray}
\end{lemma}

\begin{remark} \rm
The relations
(\ref{eq:hoq}), (\ref{eq:hoq2})
are often called the {\it higher order $q$-Serre relations} \cite[p.57]{lu}.
\end{remark}

\medskip
\noindent We now give the higher order $q$-Serre relations
in a modified form. We will use the following notation.
For all integers $n,k$ with $k$ nonnegative we define
\begin{eqnarray*}
\lbrack n\rbrack_k=
\lbrack n\rbrack
\lbrack n-1\rbrack
\cdots
\lbrack n-k+1\rbrack.
\end{eqnarray*}
We interpret $\lbrack n\rbrack_0=1$.


\begin{lemma}
\label{lem:homod}
The following relations hold in $\mathcal A$.
Let $r$ denote a positive integer
and let $n$ denote an integer at least $2r+1$.
Then for $0 \leq i\leq n$ we have
\begin{eqnarray}
\label{eq:ho}
&&x^iy^rx^{n-i}=
\sum_{\xi=0}^r
\frac{
\lbrack i\rbrack_{\xi}
\lbrack r-i\rbrack_{r-\xi}
\lbrack n-i\rbrack_{r}
}
{
\lbrack \xi \rbrack_{\xi}
\lbrack r-\xi \rbrack_{r-\xi}
\lbrack n-\xi \rbrack_{r}
}
x^{\xi}y^rx^{n-\xi}
\\
&&\qquad \qquad \qquad +
\sum_{\zeta=0}^{r-1}
\frac{
\lbrack i\rbrack_{r+1}
\lbrack i-n+r-1\rbrack_{r-\zeta-1}
\lbrack n-i\rbrack_{\zeta}
}
{
\lbrack n-\zeta \rbrack_{r+1}
\lbrack r-\zeta-1 \rbrack_{r-\zeta-1}
\lbrack \zeta \rbrack_{\zeta}
}
x^{n-\zeta}y^rx^{\zeta}.
\label{eq:zeta}
\end{eqnarray}
\end{lemma}
\noindent
{\bf Proof:}
We
abbreviate $w_i=x^iy^rx^{n-i}$
for $0 \leq i \leq n$. We obtain a system of $n-2r$ linear equations relating
$w_0,w_1,\ldots,w_n$ as follows.
For $0 \leq j\leq n-2r-1$ we multiply each term
in 
(\ref{eq:hoq}) on the left by $x^j$ and on the right by
$x^{n-2r-1-j}$. The result is a linear equation
relating
$w_j,
w_{j+1},\ldots,
w_{j+2r+1}$.
This yields a system of $n-2r$ linear equations
involving
$w_0,w_1,\ldots,w_n$.
We now solve this system 
to obtain each of
$w_0,w_1,\ldots,w_n$
in terms of $w_0,w_1,\ldots,w_r$ and
 $w_{n-r+1},w_{n-r+2},\ldots,w_n$.
We define a $(n+1)\times (2r+1)$ matrix $S$ as follows.
For
$0 \leq i \leq n$ and
$-r \leq j \leq r$ the entry
$S_{ij}=q^{2ij}$.
For $0 \leq i\leq n$ let $S_i$ denote row $i$ of $S$.
Observe that for $0 \leq j\leq n-2r-1$ the vectors
$S_j, S_{j+1}, \ldots, S_{j+2r+1}$ satisfy the same 
equation as the one for
$w_j,
w_{j+1},\ldots,
w_{j+2r+1}$ which we mentioned above.
Observe $S_0, S_1, \ldots, S_{2r}$
are linearly independent since
 $S$ is
essentially Vandermonde.
Therefore $S_0, S_1, \ldots, S_{2r}$ 
is a basis for $\K^{2r+1}$.
Consider the linear transformation $\varepsilon $
from $\K^{2r+1}$ to
$\mbox{Span}(
w_0, w_1, \ldots, w_n)$ which 
sends $S_i$ to $w_i$ for $0 \leq i \leq 2r$.
From our above comments $\varepsilon $ sends 
$S_i$ to $w_i$ for $0 \leq i \leq n$.
For $0 \leq i\leq n$ we have
\begin{eqnarray}
\label{eq:hoproof}
&&S_i=
\sum_{\xi=0}^r
\frac{
\lbrack i\rbrack_{\xi}
\lbrack r-i\rbrack_{r-\xi}
\lbrack n-i\rbrack_{r}
}
{
\lbrack \xi \rbrack_{\xi}
\lbrack r-\xi \rbrack_{r-\xi}
\lbrack n-\xi \rbrack_{r}
}
S_{\xi}
\\
&&\qquad \qquad \qquad +
\sum_{\zeta=0}^{r-1}
\frac{
\lbrack i\rbrack_{r+1}
\lbrack i-n+r-1\rbrack_{r-\zeta-1}
\lbrack n-i\rbrack_{\zeta}
}
{
\lbrack n-\zeta \rbrack_{r+1}
\lbrack r-\zeta-1 \rbrack_{r-\zeta-1}
\lbrack \zeta \rbrack_{\zeta}
}
S_{n-\zeta}.
\label{eq:zetaproof}
\end{eqnarray}
This can be verified using Cramer's rule.
Applying $\varepsilon$ to each term in 
(\ref{eq:hoproof}), 
(\ref{eq:zetaproof}) we obtain
(\ref{eq:ho}), 
(\ref{eq:zeta}).
\hfill $\Box $ \\

\begin{remark}
\label{rem:m}
 \rm
Applying each of $\sigma_0, \sigma_1,\sigma$
to the relations in Lemma
\ref{lem:homod}, we obtain additional relations involving $x, y$ which hold in 
$\mathcal A$.
\end{remark}

\noindent The following fact is an immediate consequence
of Lemma
\ref{lem:homod}
and Remark
\ref{rem:m}.

\begin{corollary}
\label{cor:dep}
The following linear dependencies hold in $\mathcal A$.
Let $r$ denote a positive integer and let $n$ denote an integer at least
$2r+1$. 
\begin{enumerate}
\item For $r+1\leq i \leq n-r$ the element $x^iy^rx^{n-i}$ is contained in
\begin{eqnarray*}
\mbox{Span}\lbrace x^{\xi}y^rx^{n-\xi} \;|\;0\leq \xi\leq r \rbrace+
\mbox{Span}\lbrace x^{n-\zeta}y^rx^{\zeta} \;|\;0\leq \zeta\leq r-1 \rbrace.
\end{eqnarray*}
\item For $r\leq i \leq n-r-1$ the element $x^iy^rx^{n-i}$ is contained in
\begin{eqnarray*}
\mbox{Span}\lbrace x^{\xi}y^rx^{n-\xi} \;|\;0\leq \xi\leq r-1 \rbrace+
\mbox{Span}\lbrace x^{n-\zeta}y^rx^{\zeta} \;|\;0\leq \zeta\leq r \rbrace.
\end{eqnarray*}
\item For $r+1\leq i \leq n-r$ the element $y^ix^ry^{n-i}$ is contained in
\begin{eqnarray*}
\mbox{Span}\lbrace y^{\xi}x^ry^{n-\xi} \;|\;0\leq \xi\leq r \rbrace+
\mbox{Span}\lbrace y^{n-\zeta}x^ry^{\zeta} \;|\;0\leq \zeta\leq r-1 \rbrace.
\end{eqnarray*}
\item For $r\leq i \leq n-r-1$ the element $y^ix^ry^{n-i}$ is contained in
\begin{eqnarray*}
\mbox{Span}\lbrace y^{\xi}x^ry^{n-\xi} \;|\;0\leq \xi\leq r-1 \rbrace+
\mbox{Span}\lbrace y^{n-\zeta}x^ry^{\zeta} \;|\;0\leq \zeta\leq r \rbrace.
\end{eqnarray*}
\end{enumerate}
\end{corollary}

\begin{definition} 
\label{def:word}
\rm
Let $n$ denote a nonnegative integer. By
a {\it word of length n} in $\mathcal A$, we mean an expression
of the form
\begin{equation}
a_1a_2\cdots a_n,
\label{eq:word}
\end{equation}
where $a_i=x$ or $a_i=y$ for $1 \leq i \leq n$.
We 
interpret the word of length 0 as the 
 identity element in $\mathcal A$.
We say this word is {\it trivial}.
By the {\it height} of the word
(\ref{eq:word}) we mean the integer
\begin{eqnarray*}
|\lbrace i\;|\;1 \leq i \leq  n, \quad a_i=x\rbrace | \;-\;
|\lbrace i\;|\;1 \leq i \leq  n, \quad a_i=y\rbrace |.
\end{eqnarray*}

\begin{definition}\rm
A word is said to be {\it balanced}
whenever it has height 0. We observe that a balanced
word has even length.
\end{definition}

\begin{example} 
We list the balanced words of length $4$.
\begin{eqnarray*}
x^2y^2,
\qquad xyxy,
\qquad xy^2x,
\qquad y^2x^2,
\qquad yxyx,
\qquad yx^2y. 
\end{eqnarray*}
\end{example}
\end{definition}

\begin{definition} \rm
Let 
$a_1a_2\cdots a_n$ 
denote a word 
in $\mathcal A$.
By the {\it height vector} of 
 $a_1a_2\cdots a_n$ 
 we mean the sequence
$(h_0,h_1,\ldots, h_n)$,
where $h_i$ denotes the height of the word
$a_{i+1}a_{i+2}\cdots a_n$ for $0 \leq i \leq n$.
We observe that the height vector of
 $a_1a_2\cdots a_n$ 
is the unique sequence  $(h_0,h_1,\ldots, h_n)$ 
 such that
(i)
$h_n=0$;
(ii)
for $1\leq i \leq n$, $h_{i-1}-h_{i}=1$
if $a_i=x$ and  $h_{i-1}-h_{i}=-1$
if $a_i=y$.
\end{definition}

\begin{example}
The word $x^2yx^3y^2x$ has height vector
\begin{eqnarray*}
(3,2,1,2,1,0,-1,0,1,0).
\end{eqnarray*}
\end{example}

\begin{definition}
\rm
Let $a_1a_2\cdots a_n$ denote a word in $\mathcal A$
and let $(h_0, h_1,\ldots, h_n)$ denote the corresponding
height vector.
We say $a_1a_2\cdots a_n$ 
is
 {\it height-symmetric}
whenever $h_i=h_{n-i}$ for $0 \leq i \leq n$.
We observe that $a_1a_2\cdots a_n$ is height-symmetric
if and only if $a_i\not=a_{n-i+1}$ for
$1\leq i \leq n$.
We remark that a height-symmetric word is balanced and invariant 
under $\sigma$.
\end{definition}

\begin{example}
We list the height-symmetric words of length 4.
\begin{eqnarray*}
x^2y^2,
\qquad xyxy,
\qquad y^2x^2,
\qquad yxyx.
\end{eqnarray*}
\end{example}

\begin{definition}\rm
\rm
Let $a_1a_2\cdots a_n$ denote a word in $\mathcal A$
and let $(h_0, h_1,\ldots, h_n)$ denote the corresponding
height vector.
We say $a_1a_2\cdots a_n$ 
is
{\it nil}
whenever at least one of 
$h_0, h_1,\ldots, h_n$ is negative.
\end{definition}

\begin{example}
We list the nil balanced words of length 4.
\begin{eqnarray*}
x^2y^2, \qquad 
xyxy, 
\qquad xy^2x,
\qquad yx^2y. 
\end{eqnarray*}
\end{example}

\medskip
\noindent We mention a few subspaces of $\mathcal A$.

\begin{definition} \rm
We let $\mathcal B$ denote the linear subspace of $\mathcal A$
spanned by the balanced words. We observe $\mathcal B$ is a $\K$-subalgebra of
$\mathcal A$ which is invariant under each of
$\sigma_0, \sigma_1, \sigma$. We let 
${\mathcal B}^{sym}$
denote the linear subspace of $\mathcal B$ spanned by the height-symmetric
words. 
We let
${\mathcal B}^{nil}$
denote the linear subspace of 
 $\mathcal B$
 spanned by the nil balanced words.
We observe 
 ${\mathcal B}^{nil}$
is a two sided ideal of $\mathcal B$ which is invariant under $\sigma$.
\end{definition}


\begin{theorem}
\label{thm:m}
Let $i,j,m,n$ denote nonnegative integers such that $i+n=j+m$. Define
$b=y^ix^my^nx^j$ and observe $b$ is a balanced word. Then
\begin{eqnarray}
b  \in 
{\mathcal B}^{sym}+
{\mathcal B}^{nil}.
\label{eq:sn}
\end{eqnarray}
\end{theorem}
\noindent {\bf Proof:}
We proceed by induction on $\mbox{min}(m,n)$. First assume
 $\mbox{min}(m,n)=0$.
Then $b$ is height-symmetric.  Therefore
$b \in {\mathcal B}^{sym}$ so
(\ref{eq:sn}) holds.
Next assume
 $\mbox{min}(m,n)>0$.
Replacing $b$ by $b^\sigma$ if necessary we may assume
$m\geq n$.
Assume for the moment that $m=n$. Then 
$b$ is height-symmetric.
Therefore
$b \in {\mathcal B}^{sym}$ so
(\ref{eq:sn}) holds.
Next assume  $m>n$. 
We may assume $n\leq j$; otherwise $b \in 
{\mathcal B}^{nil}$.
By Corollary
\ref{cor:dep}(i),
$b$ is contained in 
\begin{eqnarray}
\mbox{Span}\lbrace y^{i}x^{\xi}y^nx^{i+n-\xi} \;|\;0\leq \xi\leq n \rbrace+
\mbox{Span}\lbrace y^{i}x^{i+n-\zeta}y^{n}x^{\zeta} \;|\;0\leq \zeta\leq n-1 
\rbrace.
\label{eq:sp}
\end{eqnarray}
By the induction hypothesis, for $0 \leq \xi \leq n-1$ the word
$y^{i}x^{\xi}y^nx^{i+n-\xi}$
is contained in
${\mathcal B}^{sym}+
{\mathcal B}^{nil}$.
For $\xi=n$ the word
$y^{i}x^{\xi}y^nx^{i+n-\xi}$ is height-symmetric
and therefore contained in
${\mathcal B}^{sym}$.
For $0\leq \zeta\leq n-1 $ the word
$y^{i}x^{i+n-\zeta}y^{n}x^{\zeta}$
is contained in 
${\mathcal B}^{nil}$.
Apparently the space displayed in
(\ref{eq:sp}) is contained in
${\mathcal B}^{sym}+{\mathcal B}^{nil}$.
It follows
$b \in {\mathcal B}^{sym}+{\mathcal B}^{nil}$ as desired.
\hfill $\Box $ \\

\begin{lemma}
\label{label:sym}
For  
$b \in {\mathcal B}^{sym}+{\mathcal B}^{nil}$ we have
$b-b^{\sigma} \in {\mathcal B}^{nil}$.
\end{lemma}
\noindent {\bf Proof:}
There exists $b_0\in {\mathcal B}^{sym}$
and there exists $b_1\in {\mathcal B}^{nil}$
such that
$b=b_0+b_1$. Observe 
$b_0^{\sigma}=b_0$ so $b-b^\sigma=b_1-b_1^{\sigma}$.
By assumption $b_1 \in {\mathcal B}^{nil}$.
Recall $\sigma$ leaves ${\mathcal B}^{nil}$ invariant so
$b_1^{\sigma}\in {\mathcal B}^{nil}$.
By these comments 
$b-b^{\sigma} \in {\mathcal B}^{nil}$.
\hfill $\Box $ \\

\begin{corollary}
\label{thm:mnew}
Let $i,j,m,n$ denote nonnegative integers such that $i+n=j+m$. Define
$b=y^ix^my^nx^j$ and recall $b$ is a balanced word. Then
\begin{eqnarray*}
b - b^{\sigma} \in {\mathcal B}^{nil}.
\end{eqnarray*}
\end{corollary}
\noindent {\bf Proof:}
Immediate from Theorem
\ref{thm:m}
and 
Lemma
\ref{label:sym}.
\hfill $\Box $ \\

%
%

\begin{corollary}
\label{cor:del}
For nonnegative integers $i,j$ we have
\begin{eqnarray*}
\lbrack y^ix^i, y^jx^j \rbrack \in {\mathcal B}^{nil}.
\end{eqnarray*}
\end{corollary}
\noindent {\bf Proof:}
Apply Corollary
\ref{thm:mnew} with $m=i$ and $n=j$.
\hfill $\Box $ \\

\noindent In order to state the next theorem we make a
definition.

\begin{definition}
\rm
Let $a_1a_2\cdots a_n$ 
denote a word in $\mathcal A$.
Observe that there exists a unique sequence
 $(i_1, i_2, \ldots, i_r)$ consisting of positive
integers such that
 $a_1a_2\cdots a_n$ is one of
$x^{i_1}
y^{i_2}
x^{i_3}
\cdots
y^{i_r}$
or
$
x^{i_1}
y^{i_2}
x^{i_3}
\cdots
x^{i_r}$
or
$y^{i_1}
x^{i_2}
y^{i_3}
\cdots
x^{i_r}$
or
$
y^{i_1}
x^{i_2}
y^{i_3}
\cdots
y^{i_r}$.
We call the sequence
 $(i_1, i_2, \ldots, i_r)$ the {\it signature} of
$a_1a_2\cdots a_n$.
\end{definition}

\begin{example}
Each  of the words $yx^2y^2x$,
$xy^2x^2y$ has signature
$(1,2,2,1)$.
\end{example}

\begin{definition}
\label{def:irr}
\rm
Let $a_1a_2\cdots a_n$ denote a word in $\mathcal A$
and let $(i_1, i_2,\ldots, i_r)$
denote the corresponding signature.
We say  $a_1a_2\cdots a_n$
is {\it reducible} whenever 
there exists an integer $s$
$(2\leq s\leq r-1)$ such that
$i_{s-1}\geq i_s<i_{s+1}$.
We say a word is {\it irreducible} whenever it is
not reducible.
\end{definition}

\begin{example}
A word of length less than 4 is irreducible.
The only reducible words of length 4 are
 $xyx^2 $ and $yxy^2$. 
\end{example}

\medskip
\noindent In the following lemma we give a necessary and sufficient
condition for a given nontrivial word to be irreducible.

\begin{lemma}
\label{lem:red}
Let $a_1a_2\cdots a_n$ denote a  nontrivial word in $\mathcal A$
and let $(i_1, i_2,\ldots, i_r)$
denote the corresponding signature. 
Then the following (i), (ii)  are equivalent.
\begin{enumerate}
\item
 The word $a_1a_2\cdots a_n$
is irreducible.
\item  There exists an
integer $t$ $(1\leq t \leq r)$
such that 
\begin{eqnarray*}
i_1<i_2 < \cdots <i_{t-1}< i_t\geq i_{t+1} \geq i_{t+2}\geq \cdots \geq i_{r-1}\geq i_r.
\end{eqnarray*}
\end{enumerate}
\end{lemma}
\noindent
{\bf Proof:} 
Routine using
Definition
\ref{def:irr}.
\hfill $\Box $ \\

\medskip
\noindent 
Consider the algebra $\mathcal A$ 
as a vector space over $\K$. It turns out that
the set of irreducible words
in $\mathcal A$ forms a basis for $\mathcal A$. 
However, in order to prove
Theorem
 \ref{con:Vspanguess}
 all we need is that this set spans
 $\mathcal A$.
We will prove this much for now and use the result to obtain Theorem
 \ref{con:Vspanguess}.
For the sake of completeness, in Section 3 we will prove
that 
 the set of irreducible words
in $\mathcal A$ is a basis for $\mathcal A$.

\begin{theorem}
\label{thm:sp}
The irreducible words in $\mathcal A$ form a spanning set for $\mathcal A$.
\end{theorem}
{\bf Proof:} 
Let
$\lambda =(\lambda_1, \lambda_2, \ldots, \lambda_r)$
denote a finite sequence of positive integers.
 We say this sequence
is {\it  nonincreasing} whenever
$\lambda_{i-1}\geq \lambda_i$
for $2 \leq i \leq r$.
Let $\Psi$ denote the set consisting of the
 nonincreasing finite sequences of positive integers.
There exists a certain linear order on $\Psi$
called the 
reverse lexicographical order. This is defined as follows.
Let $\lambda =(\lambda_1, \lambda_2, \ldots, \lambda_r)$
and $\mu =(\mu_1, \mu_2, \ldots, \mu_s)$
denote elements in $\Psi$. Then
$\lambda$ is less than  $\mu$ in the reverse lexicographical order
 whenever (i) $r<s$; or
(ii) $r=s$ and there exists an integer $k$ $(1 \leq k \leq r)$
such that
$\lambda_i=\mu_i$ $(k+1\leq i \leq r)$, $\lambda_k<\mu_k$.
Let $\lambda =(\lambda_1, \lambda_2, \ldots, \lambda_r)$
denote a finite sequence of positive integers.
We let ${\overline {\lambda}}$ denote the rearrangement of
$\lambda $ into a nonincreasing sequence.
In other words
${\overline {\lambda}}= 
(\mu_1, \mu_2,\ldots, \mu_r)$ where
$\lbrace \mu_1, \mu_2,\ldots, \mu_r\rbrace =
\lbrace \lambda_1, \lambda_2, \ldots, \lambda_r\rbrace $
as  multisets and
the sequence 
$(\mu_1, \mu_2,\ldots, \mu_r)$ is nonincreasing. 
We call ${\overline {\lambda}}$
the {\it rearrangement} of $\lambda$. 
Let $\lambda =(\lambda_1, \lambda_2, \ldots, \lambda_r)$
denote a finite sequence of positive integers.
By an {\it inversion} in $\lambda $ we mean an ordered
pair $(i,j)$ of integers such that
 $1 \leq i<j\leq r$ and 
$\lambda_i<\lambda_j$. We observe $\lambda = \overline{\lambda }$
if and only if $\lambda $ has no inversions.
We assume the present theorem is false and obtain a contradiction.
By a {\it counterexample} we mean a word in $\mathcal A$
which is not contained in the span of the irreducible words.
By assumption there exists a counterexample.
Let $w$ denote a counterexample
 and let
$\lambda =(\lambda_1, \lambda_2, \ldots, \lambda_r)$ denote the corresponding 
signature.
Without loss, we assume
that among all the counterexamples, the rearrangement 
$\overline \lambda$ is minimal
with respect to the reverse lexicographical order.
Moreover, without loss we may assume
that among all the 
 counterexamples 
for which the
 rearranged signature is equal to $\overline \lambda$,
the signature
$\lambda$ has a minimal number of inversions. 
Since $w$ is a counterexample it is reducible.
Therefore there exists an integer $s$ $(2 \leq s \leq r-1)$ such 
that $\lambda_{s-1}\geq \lambda_s<\lambda_{s+1}$.
By the construction 
there exist words $w_1, w_2$ in $\mathcal A$
such that
 $w=w_1x^{\lambda_{s-1}}y^{\lambda_s}x^{\lambda_{s+1}}w_2$
or
 $w=w_1y^{\lambda_{s-1}}x^{\lambda_s}y^{\lambda_{s+1}}w_2$.
We treat the first case; the second case is treated in a similar
manner. 
By Corollary
\ref{cor:dep}(ii) (with $i=\lambda_{s-1}$, $r=\lambda_s$, 
$n=\lambda_{s-1}+\lambda_{s+1}$)
we find
$w$ is contained in 
\begin{eqnarray}
\label{eq:w1}
&&\mbox{Span}\lbrace w_1x^{\xi}y^{\lambda_s}x^{\lambda_{s-1}+\lambda_{s+1}-\xi}w_2 \;|\;0\leq \xi\leq \lambda_s-1 \rbrace\\
&&\qquad \qquad \qquad +
\quad 
\mbox{Span}\lbrace w_1x^{\lambda_{s-1}+\lambda_{s+1}-\zeta}y^{\lambda_s}x^{\zeta}w_2 \;|\;0\leq \zeta\leq \lambda_s \rbrace.
\label{x}
\end{eqnarray}
By our minimality assumptions no word which appears in 
(\ref{eq:w1}) is a counterexample
and no word which appears in
(\ref{x}) is a counterexample.
Therefore the sum displayed in
(\ref{eq:w1}),
(\ref{x}) is contained in the span of the irreducible words.
The word $w$ is contained in this space so $w$ is in the span of the
irreducible words. This contradicts our assumption that $w$ is a 
counterexample.
The result follows.
\hfill $\Box $ \\

\medskip
\noindent
We are now ready to prove 
Theorem \ref{con:Vspanguess}.

\medskip
\noindent{\bf Proof of Theorem 
\ref{con:Vspanguess}:}
By 
(\ref{eq:ms1}),
(\ref{eq:ms2}) there exists an $\mathcal A$-module structure
 on $V$
such that  $x.v=Rv$ and $y.v=Lv$ for all $v \in V$.
Let $b$ denote a balanced word in $\mathcal A$.
Using (\ref{eq:raction})
 and 
(\ref{eq:laction}) we find
each of $U_0,U_1,\ldots, U_d$ is invariant under $b$.
In particular $U_0$ is invariant under $b$.
Let $i$ denote a nonnegative integer.
The word $y^ix^i$ is balanced so it leaves $U_0$ invariant.
From 
(\ref{eq:nil})
we find $x^{d+1}$ vanishes on $V$.
In particular $x^{d+1}$ vanishes on $U_0$.
Therefore $y^ix^i$ vanishes on $U_0$ for $i>d$.
We claim the elements $y^ix^i$ $(0 \leq i \leq d)$ mutually commute on
$U_0$.
To see this, observe by
(\ref{eq:raction}),
(\ref{eq:laction}) that
each nil word of $\mathcal A$ vanishes
on $U_0$. In particular each nil word 
in
$\mathcal B$ vanishes on
$U_0$ so ${\mathcal B}^{nil}$ vanishes on $U_0$. 
By this and Corollary
\ref{cor:del}
we find $\lbrack y^ix^i,y^jx^j\rbrack$ vanishes on $U_0$
 for $0 \leq i,j\leq d$.
We have now shown the elements $y^ix^i$ $(0 \leq i \leq d)$ mutually commute
on $U_0$.
Since $\fld$ is algebraically closed there exists a nonzero 
$v \in U_0$ which is a common eigenvector of 
the $y^ix^i$ $(0 \leq i \leq d)$. 
Let $W$ denote the subspace of $V$ spanned by the
vectors
of the form
\begin{eqnarray}
\label{eq:zsp}
y^{i_1}x^{i_2}y^{i_3}x^{i_4}\cdots x^{i_n}.v,
\end{eqnarray}
where 
$i_1, i_2, \ldots, i_n$ ranges over all sequences
such that $n$ is a nonnegative even integer, and 
 $i_1, i_2, \ldots i_n$ are integers
satisfying $0\leq i_1<i_2<\cdots <i_n\leq d$.
We show $W=V$.
Observe $v \in W$ so $W\not=0$.
By Definition
\ref{def:tdp}(iv)
the space $V$ is irreducible as a module for
$A,A^*$. In order to show $W=V$ we show
$W$ is invariant under each of $A, A^*$.
Using
Lemma
\ref{lem:red}
and
 Theorem
\ref{thm:sp}
 we routinely find
${\mathcal A}.v=W$. From this we find
$W$ is invariant under each of $x,y$.
It follows $W$ is invariant under each of $R,L$.
Using
(\ref{eq:raction}),
(\ref{eq:laction})
we find the vector 
(\ref{eq:zsp})
is contained in $U_i$
where
$i=\sum_{h=1}^n i_h(-1)^h$.
By this and
(\ref{eq:fi1}),
(\ref{eq:fi2})
the vector 
(\ref{eq:zsp})
is an eigenvector for each of $F_0, F_1, \ldots, F_d$.
We now see $W$ is invariant under each of
$F_0, F_1, \ldots, F_d$.
By these comments and 
(\ref{eq:defR}), (\ref{eq:defL})
we find $W$ is invariant under each of $A, A^* $.
We conclude $W=V$ and
the result follows.
\hfill $\Box $ \\

\medskip
\noindent {\bf Proof of Theorem
\ref{thm:us}:} Combine Theorem
\ref{con:Vspanguess} with 
(\ref{eq:raction}) and
(\ref{eq:laction}).
\hfill $\Box $ \\

\medskip
\noindent {\bf Proof of Theorem 
\ref{thm:bc}:}
For $0 \leq i \leq d$, the spanning set for $U_i$ given in
Theorem
\ref{thm:us} 
has cardinality
$\bigl( {{d}\atop {i}} \bigr)$.
\hfill $\Box $ \\

\section{Comments and suggestions for further research}

\noindent
In this section we give some comments and suggestions for further research.
 We begin with a comment.

\begin{lemma}
\label{lem:fc}
The quotient algebra ${\mathcal B}/{\mathcal B}^{nil}$ is generated by
the elements
\begin{eqnarray*}
y^ix^i+{\mathcal B}^{nil} \qquad (0\leq i < \infty).
\end{eqnarray*}
Moreover 
 ${\mathcal B}/{\mathcal B}^{nil}$ is commutative.
\end{lemma}
{\bf Proof:} 
The first assertion is immediate from Theorem
\ref{thm:sp}.
The last assertion follows from this and Corollary
\ref{cor:del}.
\hfill $\Box $ \\

\medskip
\noindent Just before 
Theorem
\ref{thm:sp}
 we asserted that 
the set of irreducible words in $\mathcal A$ is a basis
for $\mathcal A$. We will now prove this assertion.

\begin{theorem}
The set of irreducible words in $\mathcal A$ is a basis
for $\mathcal A$.
\end{theorem}
{\bf Proof:}
Let $X$ denote the set of irreducible words in $\mathcal A$.
Then  $X$ spans $\mathcal A$
by Theorem
\ref{thm:sp}.
For $0 \leq n < \infty $ let 
$X_n$ denote the set of irreducible words in $\mathcal A$ which have
length $n$. 
Of course $X=\cup_{n=0}^\infty X_n$.
For $0 \leq n < \infty $ let 
${\mathcal A}_n$
denote the subspace of $\mathcal A$ spanned by the
words of length $n$.
Since the $q$-Serre relations
 are homogeneous we have
\begin{eqnarray}
\label{eq:dir}
{\mathcal A}
=\sum_{n=0}^\infty {\mathcal A}_n \qquad \qquad (\mbox{direct sum}).
\end{eqnarray}
By the construction $X_n \subseteq {\mathcal A}_n$ for 
$0 \leq n< \infty$.
By these comments we find $X_n$ spans ${\mathcal A}_n$ for
$0 \leq n < \infty $. 
We show
\begin{eqnarray}
|X_n|=\mbox{dim}({\mathcal A}_n)
\qquad \qquad 
 (0 \leq n < \infty).
\label{eq:req}
\end{eqnarray}
Let $v$ denote an indeterminate. In what follows, we consider
formal power series in $v$ which have coefficients in $\C$. 
Line 
 (\ref{eq:req}) will follow once
we  show
\begin{eqnarray}
\sum_{n=0}^\infty |X_n|v^n=\sum_{n=0}^\infty \mbox{dim}({\mathcal A}_n)v^n.
\label{eq:req2}
\end{eqnarray}
We claim 
\begin{eqnarray}
\label{eq:dimgen}
\sum_{n=0}^\infty \mbox{dim}({\mathcal A}_n)v^n=
\prod_{m=1}^\infty (1-v^{2m})^{-1}(1-v^{2m-1})^{-2}.
\end{eqnarray}
To see
(\ref{eq:dimgen}), recall $\mathcal A$ is the positive
part of 
$U_q(\widehat{sl_2})$ and hence isomorphic to the 
Verma module for 
$\widehat{sl_2}$ \cite[p. 123]{waki}
as a graded vector space.  
Apparently the left-hand side of
(\ref{eq:dimgen}) is equal to the formal character of the Verma module.
This character is inverse to the principally specialized Weyl
denominator and is therefore equal to the right-hand side of
(\ref{eq:dimgen}) \cite[p. 181]{kac}. We now have 
(\ref{eq:dimgen}).
For $0 \leq n < \infty$,
by a {\it partition of $n$} we mean a 
 sequence of positive integers
$\lambda=(\lambda_1, \lambda_2,\ldots, \lambda_r)$ 
such that $\lambda_{i-1}\geq \lambda_i$ for $2 \leq i \leq r$
and $\sum_{i=1}^r \lambda_i=n$.
We call
$\lambda_1, \lambda_2, \ldots, \lambda_r$ the {\it parts} of $\lambda$.
We sometimes write $n=|\lambda |$.
Let $p_n$ denote the number of partitions of $n$.
Let $p'_n$ denote the number of partitions of $n$ whose parts are
mutually distinct.
The following two generating functions are well known:
\begin{eqnarray*}
\sum_{n=0}^\infty p_n v^n &=&\prod_{m=1}^\infty (1-v^m)^{-1},
\\
\sum_{n=0}^\infty p'_n v^n &=&\prod_{m=1}^\infty (1-v^{2m-1})^{-1}.
\end{eqnarray*}
See for example
\cite[Theorem 1.1]{Andrews}.
 For $0 \leq n < \infty $ let $Y_n$ denote the set of ordered pairs
$(\lambda, \mu )$
 such that (i) $\lambda$ is a partition whose parts are mutually
distinct; (ii) $\mu $ is a partition; (iii) $|\lambda |+|\mu|=n$.
From the construction
\begin{eqnarray}
\sum_{n=0}^\infty |Y_n|v^n &=&
\biggl(\sum_{n=0}^\infty p'_n v^n\biggr)
\biggl(\sum_{n=0}^\infty p_n v^n\biggr)
\nonumber
\\
 &=&\prod_{m=1}^\infty 
(1-v^{2m})^{-1}
(1-v^{2m-1})^{-2}.
\label{eq:pr}
\end{eqnarray}
We show
\begin{eqnarray}
|X_n|=|Y_n| \qquad \qquad 
(0 \leq n < \infty).
\label{eq:yx}
\end{eqnarray}
To obtain 
(\ref{eq:yx}),
for $0 \leq n < \infty$
 we 
display a bijection $\varepsilon :Y_n \rightarrow X_n$. 
Let $(\lambda, \mu)$ denote an element of $Y_n$ and write
$\lambda 
=(\lambda_1, \lambda_2, \ldots, \lambda_r)$,
$\mu =(\mu_1, \mu_2, \ldots, \mu_s)$.
The image of $(\lambda, \mu)$ under $\varepsilon$ is
 the following word in $\mathcal A$:
\begin{eqnarray}
\label{eq:finword}
\cdots x^{\lambda_4}
y^{\lambda_3}
x^{\lambda_2}
y^{\lambda_1}
x^{\mu_1}
y^{\mu_2}
x^{\mu_3}
y^{\mu_4}
\cdots.
\end{eqnarray}
By the construction
the word
(\ref{eq:finword}) is irreducible with length $n$.
Therefore the word (\ref{eq:finword}) is contained in $X_n$.
Using Lemma
\ref{lem:red}
we  find $\varepsilon :
Y_n \rightarrow X_n$ is a  bijection.
We now have
(\ref{eq:yx}).
Combining
(\ref{eq:dimgen}),
(\ref{eq:pr}),
(\ref{eq:yx})
we obtain
(\ref{eq:req2}) and
(\ref{eq:req}) 
follows.
We conclude
$X_n$ is a basis of ${\mathcal A}_n$ for $0\leq n< \infty$.
By this and
(\ref{eq:dir})
we find $X$ is a basis for $\mathcal A$.
\hfill $\Box $ \\

\begin{remark}
\rm
The following is a specialization of the Jacobi triple 
product identity \cite[p. 185]{kac}:
\begin{eqnarray}
\sum_{n=-\infty}^\infty (-v)^{n^2}=
\prod_{m=1}^\infty 
(1-v^{2m})
(1-v^{2m-1})^{2}.
\label{eq:thf}
\end{eqnarray}
If we set 
 $v=e^{\pi i z}$ for $i, z \in \C$ with $i^2=-1$ and 
 $\mbox{Im}(z)>0$, then
 either side of 
(\ref{eq:thf})
is equal to the theta function $\theta_4(z)$.
See 
\cite[pp. 103--105]{Sloane} for more detail.   
\end{remark}

\medskip
\noindent We give a suggestion for future research.

\begin{problem}
\label{prob:basis}
\rm
Let $\beta, \gamma, \gamma^*, \varrho, \varrho^*$ denote
a sequence of scalars taken from $\K$.
Let $\mathcal T$ denote the associative $\K$-algebra with identity generated
by symbols $x,y$ subject to the tridiagonal relations
\begin{eqnarray*}
\lbrack x,x^2y-\beta xyx + 
yx^2 -\gamma (xy+yx)-\varrho y\rbrack &=&0, 
\\
\lbrack y,y^{2}x-\beta yxy + xy^{2} -\gamma^* (yx+xy)-
\varrho^* x\rbrack&=&0.
\end{eqnarray*}
Find analogs of the higher order $q$-Serre relations which hold in
$\mathcal T$. Find a  basis for the $\K$-vector space $\mathcal T$.
\end{problem}

\begin{problem}
\rm
Let the algebra $\mathcal T$ be as in Problem
\ref{prob:basis}. An element of $\mathcal T$
is called {\it central} whenever it commutes with every
element of $\mathcal T$. By definition 
the {\it center} of $\mathcal T$
is the $\K$-subalgebra of $\mathcal T$ consisting of
the central elements of $\mathcal T$. Describe the center
of $\mathcal T$. Find a generating set for this center.
\end{problem}

%

\section{Acknowledgment}  The authors thank 
Georgia Benkart,
Vyjayanthi Chari, 
Arjeh Cohen,
Atsushi Matsuo, Anne Schilling, and
Kenichiro Tanabe
for helpful discussions on the subject of this paper.
The authors thank Michio Jimbo for 
providing a proof of 
(\ref{eq:dimgen}).

\medskip
\noindent
Tatsuro Ito \hfil\break
Department of Computational Science\hfil\break
Faculty of Science \hfil\break
Kanazawa University \hfil\break
Kakuma-machi, Kanazawa 920--1192, Japan \hfil\break
E-mail: ito@kappa.s.kanazawa-u.ac.jp \hfil\break

\medskip
\noindent
Paul Terwilliger \hfil\break
Department of Mathematics\hfil\break
University of Wisconsin \hfil\break
480 Lincoln drive  \hfil\break
Madison, Wisconsin, 53706, USA \hfil\break
E-mail: terwilli@math.wisc.edu \hfil\break

\end{document}